\documentclass[11pt]{article}

\newcommand{\notfree}[1]{}

\newcommand{\pathtotrunk}{./}

\usepackage{amsmath,amssymb,amsfonts,amsthm}
\usepackage{ifpdf}
\usepackage{enumerate}
\usepackage{ulem}
\usepackage{leftidx}

\usepackage[affil-it]{authblk}

\ifpdf

\else

\fi

\newcommand{\arxiv}[1]{\href{http://arxiv.org/abs/#1}{\tt arXiv:\nolinkurl{#1}}}
\newcommand{\doi}[1]{\href{http://dx.doi.org/#1}{{\tt DOI:#1}}}

\newcommand{\googlebooks}[1]{(preview at \href{http://books.google.com/books?id=#1}{google books})}

\def\RCS$#1: #2 ${\expandafter\def\csname RCS#1\endcsname{#2}}
\RCS$Revision: 1913 $ \RCS$Date: 2010-07-12 09:25:35 -0700 (Mon, 12 Jul 2010) $

\theoremstyle{plain}
\newtheorem{prop}{Proposition}[section]

\newtheorem{thm}[prop]{Theorem}

\newtheorem{lem}[prop]{Lemma}

\newtheorem*{cor*}{Corollary}

\newtheorem*{defn*}{Definition}             
\newtheorem*{rem*}{Remark}               
\newtheorem*{ex*}{Example}                

\newenvironment{rem}{\noindent\textsl{Remark.}}{\\}  
\numberwithin{equation}{section}

\newcounter{comment}
\newcommand{\noop}[1]{}

\def\clap#1{\hbox to 0pt{\hss#1\hss}}

\newcommand{\Integer}{\mathbb Z}

\newcommand{\cP}{\mathcal{P}}

\newcommand{\cV}{\mathcal{V}}
\newcommand{\cW}{\mathcal{W}}

\renewcommand{\imath}{\mathfrak{i}}
\renewcommand{\jmath}{\mathfrak{j}}

\makeatletter
\newcommand{\hashdef}[2]{\@namedef{#1}{#2}}
\newcommand{\hashlookup}[1]{\@nameuse{#1}}
\makeatother

\IfFileExists{\pathtotrunk ../../graphs/lookup.tex}%
	{\newcommand{\pathtographs}{\pathtotrunk ../../graphs/}}%
	{\newcommand{\pathtographs}{diagrams/graphs/}}

\hashdef{bwd1v1duals1v1}{7BFEC24E607C77CA}%
\hashdef{bwd1v1p1p1p1duals1v4x3x2x1}{3C540D5FC4D87EC3}%
\hashdef{bwd1v1p1v1x1v1duals1v1x2v1}{434A46414835DC7C}%
\hashdef{bwd1v1v1p1p1v1x0x0p0x1x0p0x0x1duals1v1v1x2x3}{28770811E6CBBC84}%
\hashdef{bwd1v1v1p1p1v1x0x0p0x1x0p0x0x1duals1v1v2x1x3}{7BF64A4A5AFB92CF}%
\hashdef{bwd1v1v1p1p1v1x0x0p0x1x0p0x0x1duals1v1v3x2x1}{23CFC50EF547E6A4}%
\hashdef{bwd1v1v1p1p1v1x0x0p0x1x0v1x0p0x1v1x0p0x1v1x0p0x1duals1v1v2x1v2x1}{8EBBA2A373A7CC6D}%
\hashdef{bwd1v1v1p1p1v1x0x0p0x1x0v1x0p1x0duals1v1v1x2}{9FAD11930F2909AC}%
\hashdef{bwd1v1v1p1v1x0p1x0p0x1v0x1x0p0x0x1p0x1x0v1x0x0p0x1x0p0x1x0v0x0x1duals1v1v3x2x1v3x2x1}{96D75097752C76D6}%
\hashdef{bwd1v1v1p1v1x0p1x0p0x1v0x1x0p0x1x0p0x0x1v1x0x0p0x0x1p0x0x1v0x0x1duals1v1v3x2x1v3x2x1}{6BA89990A75C6ED6}%
\hashdef{bwd1v1v1p1v1x0p1x0p1x0p0x1duals1v1v4x2x3x1}{BC04D49137872F04}%
\hashdef{bwd1v1v1p1v1x0p1x0p1x0p0x1duals1v1v4x3x2x1}{E16F8C39E0489BBE}%
\hashdef{bwd1v1v1v1p1p1v0x0x1p0x0x1duals1v1v1x2x3}{CD73E7EC42BC0A68}%
\hashdef{bwd1v1v1v1p1p1v0x1x0p0x0x1v1x0p0x1duals1v1v1x2x3v2x1}{8FAB4C9171733C9F}%
\hashdef{bwd1v1v1v1p1p1v0x1x0p0x1x0duals1v1v3x2x1}{8719D37B34C923B4}%
\hashdef{bwd1v1v1v1p1p1v1x0x0p0x0x1v1x0p0x1duals1v1v3x2x1v1x2}{3C52D7BCAC7345F3}%
\hashdef{bwd1v1v1v1p1p1v1x0x0p0x0x1v1x0p0x1duals1v1v3x2x1v2x1}{DF681CAC4972E925}%
\hashdef{bwd1v1v1v1p1p1v1x0x0p0x1x0v1x0p0x1duals1v1v1x2x3v1x2}{4F97B41B2E8CAC38}%
\hashdef{bwd1v1v1v1p1p1v1x0x0p0x1x0v1x0p0x1duals1v1v1x2x3v2x1}{20FE18846E5E4AAA}%
\hashdef{bwd1v1v1v1p1p1v1x0x0p1x0x0duals1v1v1x2x3}{A52C066C3571F824}%
\hashdef{bwd1v1v1v1p1v0x1p0x1v1x0p0x1p0x1v0x0x1v1duals1v1v1x2v1x2x3v1}{AD0FD12331645520}%
\hashdef{bwd1v1v1v1p1v0x1p0x1v1x0p1x0p0x1v0x1x0v1duals1v1v1x2v1x2x3v1}{BAC8DDF76679121F}%
\hashdef{bwd1v1v1v1p1v1x0p0x1v1x0p1x1v1x0v1duals1v1v1x2v1x2v1}{48383DE8F4973A95}%
\hashdef{bwd1v1v1v1p1v1x0p0x1v1x1p0x1v0x1v1duals1v1v1x2v1x2v1}{8D4DF37334AC95F7}%
\hashdef{bwd1v1v1v1p1v1x0p0x1v1x1p1x0v0x1v1duals1v1v1x2v1x2v1}{B0A3F777BD07E47F}%
\hashdef{bwd1v1v1v1p1v1x0p1x0v1x0p1x0p0x1v1x0x0v1duals1v1v1x2v1x2x3v1}{DF170089A5590362}%
\hashdef{bwd1v1v1v1v1p1p1v1x0x0p0x1x0v1x0v1v1v1duals1v1v1v1x2v1v1}{EE553429A3976217}%
\hashdef{bwd1v1v1v1v1v1p1p1v1x0x0p0x1x0v1x0v1v1duals1v1v1v1x2x3v1v1}{5AB2F40E2972944E}%

\newcommand{\bigraph}[1]{{\hspace{-3pt}\begin{array}{c}%
  \raisebox{-2.5pt}{\includegraphics[height=6mm]{\pathtographs \hashlookup{#1}}}%
\end{array}\hspace{-3pt}}}

\newcommand{\Tr}{\operatorname{Tr}}

\ifpdf
\usepackage[pdftex,plainpages=false,hypertexnames=false,pdfpagelabels]{hyperref}
\else
\usepackage[dvips,plainpages=false,hypertexnames=false]{hyperref}
\fi

\ifpdf
\usepackage[pdftex]{graphicx}
\usepackage[pdftex]{rotating}
\else
\usepackage[dvips]{graphicx}
\usepackage[dvips]{rotating}
\fi

\usepackage{microtype}

\usepackage{tikz}
\usetikzlibrary{shapes}
\usetikzlibrary{backgrounds}

\usepackage{color}
\usepackage{xcolor}
\definecolor{dark-red}{rgb}{0.7,0.25,0.25}
\definecolor{dark-blue}{rgb}{0.15,0.15,0.55}
\definecolor{medium-blue}{rgb}{0,0,0.65}
\hypersetup{
   colorlinks, linkcolor={purple},
   citecolor={medium-blue}, urlcolor={medium-blue}
}

\title{Subfactors of index exactly 5}

\author{Masaki Izumi}
\affil{Kyoto University}

\author{Scott Morrison}
\affil{Mathematical Sciences Institute, Australian National University}

\author{David Penneys}
\affil{University of Toronto}

\author{Emily Peters}
\affil{Loyola University Chicago}

\author{Noah Snyder}
\affil{Indiana University}

\usepackage{fullpage}

\begin{document}

\maketitle

\begin{abstract}
We give the classification of subfactor planar algebras at index exactly 5. 
All the examples arise as standard invariants of subgroup subfactors. 
Some of the requisite uniqueness results come from work of Izumi in preparation. 
The non-existence results build upon the classification of subfactor planar algebras with index less than 5, with some additional analysis of special cases.
\end{abstract}


\section{Introduction}

The goal of this article is to extend our prior classification of subfactors of index strictly less than $5$ \cite{MR2914056,MR2902285,MR2993924,MR2902286} to include subfactors of index $5$.  
Recall that any finite index subfactor $N \subset M$ has a `standard invariant' which plays the role of a Galois group.  
For an arbitrary factor $M$ it may be difficult to answer the ``inverse Galois problem" of determining which standard invariants are realized and in how many ways, but by deep results of Ocneanu and Popa \cite{MR996454,MR1055708,MR1278111}, when $M$ is the hyperfinite {\rm II}$_1$ factor and the standard invariant is finite depth or amenable, the standard invariant is a complete invariant.  
In this article we classify all subfactor standard invariants with index exactly $5$.

\begin{thm}
\label{MainTheorem}
Besides subfactors with standard invariant Temperley-Lieb (and principal graph $A_\infty$) and the non-extremal perturbation of the $A_\infty^{(1)}$ subfactor, there are seven subfactor standard invariants at index exactly 5.
There are five subgroup subfactors
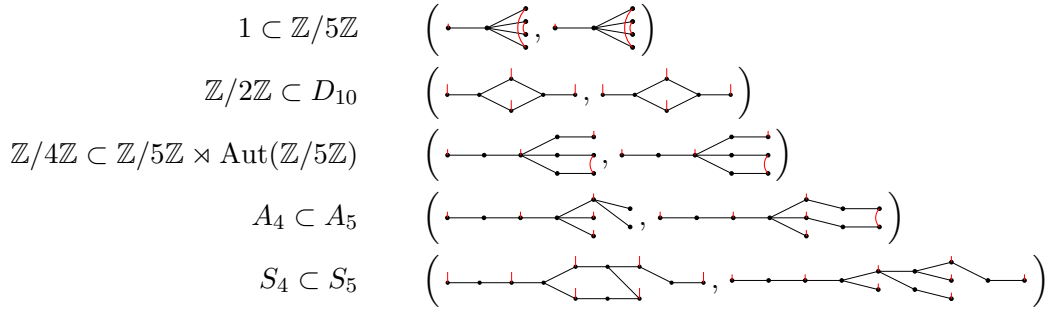
\begin{align*}
1 \subset \Integer/5\Integer & \qquad \left(\bigraph{bwd1v1p1p1p1duals1v4x3x2x1}, \bigraph{bwd1v1p1p1p1duals1v4x3x2x1}\right) \\
\Integer/2\Integer \subset D_{10} & \qquad \left(\bigraph{bwd1v1p1v1x1v1duals1v1x2v1}, \bigraph{bwd1v1p1v1x1v1duals1v1x2v1}\right) \\
\Integer/4\Integer \subset \Integer/5\Integer \rtimes \operatorname{Aut}(\Integer/5\Integer) & \qquad \left(\bigraph{bwd1v1v1p1p1v1x0x0p0x1x0p0x0x1duals1v1v2x1x3}, \bigraph{bwd1v1v1p1p1v1x0x0p0x1x0p0x0x1duals1v1v2x1x3}\right) \\
A_4 \subset A_5 & \qquad \left(\bigraph{bwd1v1v1v1p1p1v0x0x1p0x0x1duals1v1v1x2x3}, \bigraph{bwd1v1v1v1p1p1v0x1x0p0x0x1v1x0p0x1duals1v1v1x2x3v2x1}\right) \\
S_4\subset S_5 & \qquad \left(\bigraph{bwd1v1v1v1p1v1x0p0x1v1x1p0x1v0x1v1duals1v1v1x2v1x2v1},\bigraph{bwd1v1v1v1p1v0x1p0x1v1x0p0x1p0x1v0x0x1v1duals1v1v1x2v1x2x3v1}\right)
\end{align*}
along with the (non-isomorphic) duals of $A_4 \subset A_5$ and $S_4 \subset S_5$.
\end{thm}

This argument follows a similar outline as previous small index classification theorems as pioneered by Jones and Ocneanu for subfactors of index less than $4$ \cite{MR0696688,JonesICM,MR996454,MR999799}. 
That classification begins by enumerating all possible principal graphs of index less than $4$ (the ADE Coxeter-Dynkin diagrams).
Then for the diagrams which don't come from subfactors ($D_{\text{odd}}$ and $E_7$), one finds obstructions that prove that no such subfactor can exist \cite{MR996454}.
Finally, for the remaining diagrams, one constructs associated subfactors and proves a uniqueness result \cite{MR0696688,MR996454,MR999799, MR1193933, MR1145672, MR1313457, MR1308617}.
 
Above index $4$, we cannot separate the enumeration and obstruction steps.
As developed by Haagerup \cite{MR1317352}, one can use certain kinds of obstructions (Ocneanu's triple point obstruction and associativity) during the enumeration step to again get a manageable list.  Haagerup's list was used to give a complete classification up to index $3+\sqrt{3}$ in \cite{MR1317352, MR1625762, MR1686551, MR2472028, 0909.4099}.  Haagerup's enumeration technique was further improved and automated in \cite{MR2914056} to give a new list up to (but not including) index $5$, and in \cite{MR2902285,MR2993924,MR2902286} we turned this into a complete classification.
Throughout this article, we assume familiarity with the techniques of these earlier classification articles, and we recommend \cite{MR3166042} for a general overview.

In the present article, the existence results are immediate, since all the principal graphs in the main theorem are realized by subgroup subfactors.  
Uniqueness of each of these realizations is proved in Section \ref{sec:uniqueness} using a mix of connection techniques and techniques from \cite{IzumiGoldman}.  
The initial enumeration of candidate principal graphs is essentially identical to \cite{MR2914056} as explained in the beginning of Section \ref{sec:nonexistence}.   
In the rest of Section \ref{sec:nonexistence}, we eliminate the remaining candidate graph pairs.
In addition to the techniques of the  classification below index 5 (as improved in \cite{MR3157990,1307.5890}), we use results of Izumi on $2^n$ and $3^n1$ star graphs \cite{IzumiGoldman} and a result of Schou describing which 4-star graphs $\Gamma$ have no connections on $(\Gamma,\Gamma)$ \cite{1304.5907}.

\subsection*{Acknowledgements}
Scott Morrison would like to thank the Erwin Schr\"odinger Institute and its 2014 programme on ``Modern Trends in Topological Quantum Field Theory'' for their hospitality.
All the authors would like to thank the Banff International Research Station for hosting the 2014 workshop on Subfactors and Fusion Categories.

Masaki Izumi was supported by JSPS, the Grant-in-Aid for Scientific Research (B) 22340032.
Scott Morrison was supported by an Australian Research Council Discovery Early Career Researcher Award, DE120100232 and Discovery Project `Subfactors and symmetries' DP140100732.
David Penneys was supported in part by the Natural Sciences and Engineering Research Council of Canada.
The last 4 authors were supported by DOD-DARPA grant HR0011-12-1-0009.


\section{Uniqueness of examples at index 5} \label{sec:uniqueness}

\subsection{Izumi's Goldman-type theorems}

We begin with two lemmas due to Izumi which will appear in \cite{IzumiGoldman} which completely classify subfactors with principal graphs $2^n$ and $3^n1$.

\begin{thm}\label{thm:2n}
If one of the principal graphs of a subfactor $A\subset B$ is the $2^n$ spoke graph
$$
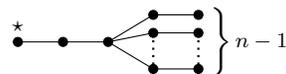
\begin{tikzpicture}[baseline=-.1cm,scale=.6]
	\filldraw (1,0) circle (.1cm) node [above] {{\scriptsize{$\star$}}};
	\filldraw (2,0) circle (.1cm);
	\filldraw (3,0) circle (.1cm) ;
	\filldraw (4,.6) circle (.1cm);
	\filldraw (4,.2) circle (.1cm);
	\filldraw (4,-.6) circle (.1cm);
	\filldraw (5,.6) circle (.1cm);
	\filldraw (5,.2) circle (.1cm);
	\filldraw (5,-.6) circle (.1cm);
	\draw[] (1,0)--(3,0);
	\draw (5,.6)--(4,.6)--(3,0)--(4,.2)--(5,.2);
	\draw (5,-.6)--(4,-.6)--(3,0);
	\draw[thick, dotted]  (4,.2)--(4,-.6);
	\draw[thick, dotted]  (5,.2)--(5,-.6);
	\node at (5.5,0) {$\bigg\}$};
	\node at (6.4,0) {\scriptsize{$n-1$}};
\end{tikzpicture}
$$
then $n=q-1$ with $q$ a prime power, and $A\subset B$ is the subgroup subfactor for $\mathbb{F}_q^{\times} \subset \mathbb{F}_q\rtimes \mathbb{F}_q^{\times}$, where  $\mathbb{F}_q$ is a finite field of order $q$.
\end{thm}

\begin{rem}
The action of $\mathbb{F}_q\rtimes \mathbb{F}_q^{\times}$ on $\mathbb{F}_q$  as an affine transformation group is a typical example of a sharply 2-transitive permutation group action.
Izumi's article \cite{IzumiGoldman} will contain a Goldman-type result for general sharply 2-transitive permutation group actions.
\end{rem}

For $3^n1$, the subfactor is not self-dual and we need a more precise statement.

\begin{thm}\label{thm:3n1}
If the graph between the $A-A$ and $B-A$ bimodules is the $3^n1$ spoke graph
$$
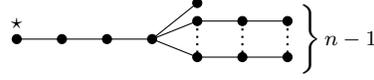
\begin{tikzpicture}[baseline=-.1cm,scale=.6]
	\filldraw (0,0) circle (.1cm) node [above] {{\scriptsize{$\star$}}};
	\filldraw (1,0) circle (.1cm);
	\filldraw (2,0) circle (.1cm);
	\filldraw (3,0) circle (.1cm) ;
	\filldraw (4,.8) circle (.1cm);
	\filldraw (4,.4) circle (.1cm);
	\filldraw (4,-.4) circle (.1cm);
	\filldraw (5,.4) circle (.1cm);
	\filldraw (5,-.4) circle (.1cm);
	\filldraw (6,.4) circle (.1cm);
	\filldraw (6,-.4) circle (.1cm);
	\draw[] (0,0)--(3,0);
	\draw (4,.8)--(3,0)--(4,.4)--(6,.4);
	\draw (6,-.4)--(4,-.4)--(3,0);
	\draw[thick, dotted]  (4,.4)--(4,-.4);
	\draw[thick, dotted]  (5,.4)--(5,-.4);
	\draw[thick, dotted]  (6,.4)--(6,-.4);
	\node at (6.5,0) {$\bigg\}$};
	\node at (7.4,0) {\scriptsize{$n-1$}};
\end{tikzpicture}
$$
then $n=q-1$ with a prime power $q$, and the subfactor is the subgroup subfactor for $F_q\rtimes F_q^{\times} \subset PGL(2,q)$.
\end{thm}

\begin{rem}
Note that $A_4\subset A_5$ is the case of $q=4$.
$PGL(2,q)$ acting on the projective geometry $PG_2(q)$ is a typical example of a sharply 3-transitive permutation group action.
Izumi's article \cite{IzumiGoldman} will contain a Goldman-type result for general sharply 3-transitive permutation group actions.
\end{rem}

\subsection{Examples at index 5}
We now give arguments for uniqueness for the examples at index 5.

\begin{prop}
These 2 principal graphs are realised by unique subfactors, namely the subgroup subfactors from Theorem \ref{MainTheorem}.
\begin{align*}
1 \subset \Integer/5\Integer & \qquad \left(\bigraph{bwd1v1p1p1p1duals1v4x3x2x1}, \bigraph{bwd1v1p1p1p1duals1v4x3x2x1}\right) \\
\Integer/2\Integer \subset D_{10} & \qquad \left(\bigraph{bwd1v1p1v1x1v1duals1v1x2v1}, \bigraph{bwd1v1p1v1x1v1duals1v1x2v1}\right) 
\end{align*}
\end{prop}
\begin{proof}
This follows from Izumi's result in \cite{MR1491121}; we additionally mention alternate proofs.
The first is depth 2 and simply laced, so it must come from a group subfactor \cite{MR860811,MR983334,MR1111570}.
The second is the graph called $D_{5}/\Integer_2$ in \cite{MR1418507}, which shows it must come from the claimed subgroup subfactor.
We also remark that the second is uniquely realized via the classification of singly generated subfactor planar algebras with $\dim(\cP_{3,\pm})\leq 13$ \cite{MR1972635}.
\end{proof}

\begin{prop}
These 2 principal graphs are realised by unique subfactors, namely the subgroup subfactors from Theorem \ref{MainTheorem}.
\begin{align*}
\Integer/4\Integer \subset \Integer/5\Integer \rtimes \operatorname{Aut}(\Integer/5\Integer) & \qquad \left(\bigraph{bwd1v1v1p1p1v1x0x0p0x1x0p0x0x1duals1v1v2x1x3}, \bigraph{bwd1v1v1p1p1v1x0x0p0x1x0p0x0x1duals1v1v2x1x3}\right) \\
A_4 \subset A_5 & \qquad \left(\bigraph{bwd1v1v1v1p1p1v0x0x1p0x0x1duals1v1v1x2x3}, \bigraph{bwd1v1v1v1p1p1v0x1x0p0x0x1v1x0p0x1duals1v1v1x2x3v2x1}\right) 
\end{align*}
\end{prop}
\begin{proof}
The first is as claimed by Theorem \ref{thm:2n}, and the second is as claimed by Theorem \ref{thm:3n1}.
\end{proof}

\begin{rem} 
It is possible to lift the shading on the 2222 principal graph pair to get a $\mathbb{Z}/2\mathbb{Z}$-graded unitary fusion category.
One can show this as follows.

There is a 1-parameter family of bi-unitary connections on the 2222 principal graph pair up to gauge equivalence and graph automorphism.
Each such connection is equivalent to one where the diagrammatic branch matrix \cite[Section 5.2]{MR2993924} is given by
$$
U=
\frac{1}{4}
\left(
\begin{array}{cccc}
 -1 & \sqrt{5} & \sqrt{5} & \sqrt{5}\\
\sqrt{5} & -\sqrt{5}\alpha &\eta \alpha & \sqrt{5} \eta \\
\sqrt{5}& 1 & \sqrt{5} \eta & -\sqrt{5} \eta \\
\sqrt{5} & \sqrt{5} \alpha & -\sqrt{5}\alpha & 1
\end{array}
\right)
$$
where $\eta\in \mathbb{T}$ is the parameter, and $\alpha=\left(\sqrt{5} \eta-1\right)/\left(\sqrt{5}-\eta\right)$ (which has modulus 1).
However, by Ocneanu rigidity \cite{MR2183279}, not all of these connections are flat.

As in \cite{MR2993924}, we look at $\Tr(UU^T)-2$, which must be equal to the sum of the eigenvalues for the new low-weight rotational eigenvectors $A,B$ at depth 3.
A straightforward calculation shows that there are only solutions for $\eta\in \mathbb{T}$ when $\Tr(UU^T)-2\in \{-1,2\}$.
When $\omega_A=\omega_B=1$, there are no flat generators at depth 3 in the graph planar algebra.
When $\omega_A=\overline{\omega_B}=\exp(2\pi i/ 3)$, there are two corresponding flat generators $A,B$.
Using the technique of \cite{1208.3637}, we have shown that $A,B$ generate a symmetrically self-dual, self-conjugate subfactor planar algebra.
Hence we may lift the shading to obtain a $\Integer/2\Integer$-graded unitary fusion category.
This adds to the limited supply of known weakly integral, non-integral fusion categories.
\end{rem}

\begin{prop}
The principal graphs below are realised by a unique subfactor, namely $S_4 \subset S_5$, because there is a unique bi-unitary connection on the principal graphs up to gauge equivalence.
\begin{align*}
S_4\subset S_5 & \qquad \left(\bigraph{bwd1v1v1v1p1v1x0p0x1v1x0p1x1v1x0v1duals1v1v1x2v1x2v1}, \bigraph{bwd1v1v1v1p1v0x1p0x1v1x0p1x0p0x1v0x1x0v1duals1v1v1x2v1x2x3v1}\right)
\end{align*}
\end{prop}
\begin{proof}
In solving the bi-unitary conditions, we see that all the equations are encoded as the unitarity of 3-by-3 or smaller matrices. One easily determines the norms of all the entries of the connection using the renormalization axiom.
When the highest valence in a principal graph is 3, one can then solve all the bi-unitary conditions using the cosine rule to express some phases in terms of others, producing only finitely many possibilities. In this case, we find only one possibility. It remains to determine how these solutions fall in gauge group orbits, which we do by an explicit calculation, exhibiting a gauge group element bringing any solution to a chosen one. The detailed calculation is available as a Mathematica notebook {\tt S4S5.nb} with the {\tt arXiv} sources of this article.
\end{proof}

\begin{rem}
Izumi also has Goldman-type results for the subgroup subfactors $S_4\subset S_5$ and $A_5\subset A_6$.
\end{rem}

\section{Enumerating graph pairs} \label{sec:nonexistence}

\subsection{Notation}

The previous classification of subfactors with index less than 5 used the terminology of weeds and vines to describe infinite families of graph pairs obtained from enumerating graph pairs satisfying Ocneanu's triple point obstruction and the associativity test \cite{MR2914056}.
A \underline{vine} is a graph pair which represents an infinite family of finite graph pairs obtained by translation.
A \underline{translation} of a graph pair means increasing the supertransitivity by an even amount. 
A \underline{weed}  is a graph pair which represents an infinite family of finite graph pairs obtained by translation and extension.
An \underline{extension} of a graph pair adds new vertices and edges at strictly greater depths than the maximum depth of any vertex in the original pair. 

Recall from \cite{MR1334479} that a principal graph pair is \underline{stable} at depth $n$ if every vertex at depth $n$ connects to at most one vertex at depth $n+1$ by at most one edge, and every vertex at depth $n+1$ is connected to exactly one vertex at depth $n$.
By \cite{MR1334479,MR3157990}, when $\delta>2$, if the principal graphs $(\Gamma_+,\Gamma_-)$ are both stable at depth $n$, or $\Gamma_+$ is stable at depths $n$ and $n+1$, then both $\Gamma_\pm$ are stable at depth $k$ for all $k\geq n+1$, and both $\Gamma_\pm$ are finite.

Recall from \cite{1310.8566} that a \underline{cylinder} is a graph pair which represents an infinite family of (finite!) graph pairs obtained by translation and stable extension, i.e., an extension which is stable at all higher depths.

\subsection{Enumerating graph pairs}

We now read the proofs of \cite{MR2914056}, checking for  places where we used the fact that the index was strictly less than 5. 
We obtain the following.

\begin{prop}\label{prop:PossibleGraphs}
Subfactors with index in the interval $(4,5]$ are either
\begin{itemize}
\item 1-supertransitive,
\item an extension of one of the  vines $\cV_\infty$  of Theorem 6.1 of \cite{MR2914056},
\item a translated extension of one of the weeds $\cW_\infty$ of Theorem 6.1 of \cite{MR2914056}, 
\item a translated extension (the extension must be stable, by \cite{MR1334479,MR3157990}) ) of the cylinder 
\begin{align*}
\Gamma_{4621} & = \left(\bigraph{bwd1v1v1v1v1p1p1v1x0x0p0x1x0v1x0v1v1v1duals1v1v1v1x2v1v1}, \bigraph{bwd1v1v1v1v1p1p1v1x0x0p0x1x0v1x0v1v1v1duals1v1v1v1x2v1v1}\right),
\end{align*}
\item
 or have one of the following principal graphs (with index exactly 5):
\begin{align*}
  G_1 & = \left(\bigraph{bwd1v1v1p1v1x0p1x0p1x0p0x1duals1v1v4x2x3x1}, \bigraph{bwd1v1v1p1v1x0p1x0p1x0p0x1duals1v1v4x2x3x1}\right) & \text{(from \S 6.2),}\displaybreak[1]\\
    G_2 & =  \left(\bigraph{bwd1v1v1p1v1x0p1x0p1x0p0x1duals1v1v4x3x2x1}, \bigraph{bwd1v1v1p1v1x0p1x0p1x0p0x1duals1v1v4x2x3x1}\right) & \text{(from \S 6.2),}\displaybreak[1]\\
    G_3 & =  \left(\bigraph{bwd1v1v1p1v1x0p1x0p1x0p0x1duals1v1v4x3x2x1}, \bigraph{bwd1v1v1p1v1x0p1x0p1x0p0x1duals1v1v4x3x2x1}\right) & \text{(from \S 6.2),} \\
     G_4 & =   \left(\bigraph{bwd1v1v1p1v1x0p1x0p0x1v0x1x0p0x1x0p0x0x1v1x0x0p0x0x1p0x0x1v0x0x1duals1v1v3x2x1v3x2x1}, \bigraph{bwd1v1v1p1v1x0p1x0p0x1v0x1x0p0x0x1p0x1x0v1x0x0p0x1x0p0x1x0v0x0x1duals1v1v3x2x1v3x2x1}\right) & \text{(2nd graph from $\cV_{11,c}$),} \displaybreak[1]\\
     G_5 & =  \left(\bigraph{bwd1v1v1v1p1v1x0p1x0v1x0p1x0p0x1v1x0x0v1duals1v1v1x2v1x2x3v1}, \bigraph{bwd1v1v1v1p1v1x0p0x1v1x1p1x0v0x1v1duals1v1v1x2v1x2v1}\right) & \text{(3rd last graph from $\cV_{10}$),} \\
      G_6 & =  \left(\bigraph{bwd1v1v1v1p1p1v1x0x0p1x0x0duals1v1v1x2x3}, \bigraph{bwd1v1v1v1p1p1v1x0x0p1x0x0duals1v1v1x2x3}\right)& \text{(from $\cV_{o2,a}$),} \displaybreak[1]\\
 G_7 & =  \left(\bigraph{bwd1v1v1v1p1p1v1x0x0p0x1x0v1x0p0x1duals1v1v1x2x3v1x2}, \bigraph{bwd1v1v1v1p1p1v1x0x0p0x1x0v1x0p0x1duals1v1v1x2x3v1x2}\right)& \text{(from $\cV_{o2,a}$),} \displaybreak[1]\\
 G_8 & =  \left(\bigraph{bwd1v1v1v1p1p1v1x0x0p0x1x0v1x0p0x1duals1v1v1x2x3v2x1}, \bigraph{bwd1v1v1v1p1p1v1x0x0p1x0x0duals1v1v1x2x3}\right)& \text{(from $\cV_{o2,a}$),}\displaybreak[1]\\
  G_9 & =  \left(\bigraph{bwd1v1v1v1p1p1v0x1x0p0x1x0duals1v1v3x2x1}, \bigraph{bwd1v1v1v1p1p1v0x1x0p0x1x0duals1v1v3x2x1}\right)& \text{(from $\cV_{o2,c}$),} \displaybreak[1]\\
 G_{10} & =  \left(\bigraph{bwd1v1v1v1p1p1v1x0x0p0x0x1v1x0p0x1duals1v1v3x2x1v1x2}, \bigraph{bwd1v1v1v1p1p1v0x1x0p0x1x0duals1v1v3x2x1}\right)& \text{(from $\cV_{o2,c}$),} \displaybreak[1]\\
 G_{11} & =  \left(\bigraph{bwd1v1v1v1p1p1v1x0x0p0x0x1v1x0p0x1duals1v1v3x2x1v2x1}, \bigraph{bwd1v1v1v1p1p1v1x0x0p0x0x1v1x0p0x1duals1v1v3x2x1v2x1}\right)& \text{(from $\cV_{o2,c}$),}\displaybreak[1]\\
  G_{12} & =  \left(\bigraph{bwd1v1v1p1p1v1x0x0p0x1x0p0x0x1duals1v1v1x2x3}, \bigraph{bwd1v1v1p1p1v1x0x0p0x1x0p0x0x1duals1v1v1x2x3}\right)& \text{(from $\cV_{e2}$),} \displaybreak[1]\\
 G_{13} & =  \left(\bigraph{bwd1v1v1p1p1v1x0x0p0x1x0p0x0x1duals1v1v3x2x1}, \bigraph{bwd1v1v1p1p1v1x0x0p0x1x0p0x0x1duals1v1v3x2x1}\right)& \text{(from $\cV_{e2}$),} \displaybreak[1]\\
 G_{14} & =  \left(\bigraph{bwd1v1v1p1p1v1x0x0p0x1x0v1x0p1x0duals1v1v1x2}, \bigraph{bwd1v1v1p1p1v1x0x0p0x1x0v1x0p1x0duals1v1v1x2}\right)& \text{(from $\cV_{e2}$),} \displaybreak[1]\\
 G_{15} & =  \left(\bigraph{bwd1v1v1p1p1v1x0x0p0x1x0v1x0p0x1v1x0p0x1v1x0p0x1duals1v1v2x1v2x1}, \bigraph{bwd1v1v1p1p1v1x0x0p0x1x0v1x0p0x1v1x0p0x1v1x0p0x1duals1v1v2x1v2x1}\right)& \text{(from $\cV_{e2}$),}\displaybreak[1]\\
  \Gamma_{5521} & =  \left(\bigraph{bwd1v1v1v1v1v1p1p1v1x0x0p0x1x0v1x0v1v1duals1v1v1v1x2x3v1v1}, \bigraph{bwd1v1v1v1v1v1p1p1v1x0x0p0x1x0v1x0v1v1duals1v1v1v1x2x3v1v1}\right) & \text{(from Theorem 6.9).}
\end{align*}
\end{itemize}
\end{prop}
\begin{proof}
First, note that although the `classification statements' described in \cite{MR2914056} purport to restrict the principal graphs of subfactors with index strictly less than some limit $\Lambda$, in fact all the odometer results remain true for describing principal graphs with index no greater than $\Lambda$.

In \cite{MR2914056} we used a special argument in Section 5 to eliminate all 1-supertransitive subfactors with index strictly less than 5. 
That argument breaks at index 5, but will be replaced with a stronger argument below in Proposition \ref{prop:1ST}.

Thus the proof of \cite[Lemma 6.2]{MR2914056} proves 
that all principal graphs of 2-supertransitive subfactors with index at most 5 are translated extensions of one of the weeds in $\cW_1$ (which is defined in that lemma).
Note \cite[Lemma 6.2]{MR2914056} is already stated as a classification statement for index at most 5, and  \cite[Lemma 6.4]{MR2914056} is not concerned with indices.

Note, now, the remark at the end of \cite[Section 6.2]{MR2914056}. 
There we discarded 3 weeds with index exactly 5, which we now have to take into account, as the first three exceptions listed in the statement above.

Now \cite[Sections 6.2.1-6.2.3]{MR2914056} already allowed for the index being exactly 5. 
However, we have to be a little careful as in \cite[Sections 6.2.2]{MR2914056} the second graph in $\cV_{11,c}$ had index 5, and similarly in \cite[Sections 6.2.3]{MR2914056} the third last graph in $\cV_{10}$ has index 5. 
As these were not carried over to the collection $\cV_\infty$ of vines in \cite[Theorem 6.1]{MR2914056}, we need to include them as the fourth and fifth exceptions listed above.

At this point all principal graphs which have supertransitivity greater than 1 and start with a triple point have been accounted for. 
We now begin reading \cite[Section 6.3]{MR2914056}.

Note \cite[Theorems 6.6-6.9]{MR2914056} already allow for the index being exactly 5; however, in each case there are some vines which we list above. 
Reading the proof of \cite[Theorem 6.10]{MR2914056}, we see that it does not rule out an index 5 principal graph coming from the weed $\Gamma_{4621}$.
Finally, \cite[Theorem 6.11]{MR2914056} already allows for the index being exactly 5, and results in no extra cases.
\end{proof}

Our task is to whittle down the possibilities arising in Proposition \ref{prop:PossibleGraphs}, until all that remains are the graphs from Theorem \ref{MainTheorem}.
We start with the 1-supertransitive possibilities, and we then consider the vines $\cV_\infty$, the weeds $\cW_\infty$, and the additional cylinder.

\begin{prop}
\label{prop:1ST}
The only possible principal graphs for 1-supertransitive subfactors with index exactly 5 are:
\begin{align*}
1 \subset \Integer/5\Integer & \qquad \left(\bigraph{bwd1v1p1p1p1duals1v4x3x2x1}, \bigraph{bwd1v1p1p1p1duals1v4x3x2x1}\right) \\
\Integer/2\Integer \subset D_{10} & \qquad \left(\bigraph{bwd1v1p1v1x1v1duals1v1x2v1}, \bigraph{bwd1v1p1v1x1v1duals1v1x2v1}\right)
\end{align*}
\end{prop}
\begin{proof}
By \cite[Theorem 5.1]{MR2914056}, either all the vertices at depth 2 have dimension 1, and thus the graph is the depth 2 graph corresponding to $1\subset \Integer/5\Integer$, or there are 2 vertices at depth 2 with dimension 2.
In the latter case, each vertex must connect to a vertex of dimension $\sqrt{5}$.
If these vertices are distinct, then the graph has annular multiplicities $*10$, and thus the dual graph has a univalent vertex at depth 1 \cite{MR2972458}.
This vertex gives a normalizer which yields an intermediate subfactor, a contradiction, since the index is not composite.
Hence both vertices at depth 2 connect to a single vertex at depth 3 with dimension $\sqrt{5}$, which must connect to a single vertex at depth 4 of dimension 1.
In this case, we must have $\dim(\cP_{3,\pm})=13$, which were completely classified by \cite{MR1972635}.
There is a unique subfactor planar algebra with this principal graph, with the dual data given above.
\end{proof}

Following \cite{MR2786219}, the index of a finite depth subfactor is a cyclotomic integer \cite{MR1120140,MR1266785,MR2183279}, which reduces each vine to finitely many possible principal graphs. 

\begin{prop}[{\cite[Section 3.3]{MR2902286}}]
The only possibilities at index 5 coming from the vines in $\cV_\infty$ are $\Gamma_{5521}$ and
\begin{align*}
G_8 = A_4\subset A_5 & =  \left(\bigraph{bwd1v1v1v1p1p1v1x0x0p0x1x0v1x0p0x1duals1v1v1x2x3v2x1}, \bigraph{bwd1v1v1v1p1p1v1x0x0p1x0x0duals1v1v1x2x3}\right).
\end{align*}
\end{prop}

\begin{prop}
The weeds in $\cW_\infty$ have already been ruled out, at all indices, in \cite{MR2902285, MR2993924}. 
Recall that the only subfactors represented by these weeds are the GHJ 3311 subfactors at index $3+\sqrt{3}$.
\end{prop}

Recall that a 4-star is a simply laced graph with a single vertex of valence 4 and all other vertices having valence at most 2.
We use the notation $S(a,b,c,d)$ for a 4-star with arms having $a,b,c,d$ edges respectively emanating from the 4-valent vertex.
If a 4-star is the principal graph of a subfactor, then the $\star$'d vertex must be on the longest arm, e.g., see \cite[Proposition 1.17]{1310.8566}.

\begin{thm}[{\cite[p. 41]{1304.5907}}]
Below is a complete list of 4-stars $\Gamma=S(a,b,c,d)$ for which the pair $(\Gamma,\Gamma)$ has a biunitary connection, which is a necessary condition for the existence of a subfactor with principal graphs $(\Gamma,\Gamma)$. 
\begin{itemize}
\item
$S(j,j,k,k)$ for $1\leq j\leq k$
\item
$S(j,j+1,j+1,j+m)$ for $1\leq j$ and $1\leq m\leq 3$
\item
$S(1,2,2,5)$
\item
$S(j,j+1,j+2,j+m)$ for $1\leq j$ and $2\leq m\leq 4$
\item
$S(j,j+2,j+2,j+2)$ for $1\leq j$
\end{itemize}
\end{thm}

\begin{prop}
There are no subfactors which are translated extensions of the cylinder $\Gamma_{4621}$.
Moreover, $\Gamma_{5521}$ is not the principal graphs of a subfactor.
\end{prop}
\begin{proof}
The graph $\Gamma_{5521}$ and all stable translated extensions of $\Gamma_{4621}$ are 4-stars of the form $S(a,b,2,1)$ where $a,b>2$.
By Schou's thesis, the only such 4-stars $\Gamma$ with bi-unitary connections on $(\Gamma,\Gamma)$ are of the form $S(j,j+1,j+2,j+m)$ for $j=1$ and $2\leq m\leq 4$, yielding $S(1,2,3,n)$ for $3\leq n\leq 5$.
\end{proof}

\subsection{Eliminating non-principal graphs}

In this subsection, we eliminate all graphs from Proposition \ref{prop:PossibleGraphs} which do not appear in Theorem \ref{MainTheorem}.
That is, we must eliminate all the $G_i$'s except $G_5,G_8,G_{13}$.

\begin{lem}
The graphs $G_{12}$, are not the principal graphs of a subfactor.
\end{lem}
\begin{proof}
This follows from Theorem \ref{thm:2n}.
\end{proof}

\begin{lem}
The graphs $G_1, G_2$ and $G_3$ are not principal graphs of subfactors.
\end{lem}
\begin{proof}
In each graph, there is a dual pair of depth 4 objects with differing dimensions (1 and 2).
\end{proof}

\begin{lem}
The graphs $G_4, G_{14}$ and $G_{15}$ are not the principal graph of a subfactor.
\end{lem}
\begin{proof}
There are vertices with dimensions less than 1.
\end{proof}

\begin{lem}
The graphs $G_7 $ and $G_{10}$ are not principal graphs of subfactors.
\end{lem}
\begin{proof}
Each graph has three vertices of dimension 1. These must form a group under tensor product, but each is self-dual, and there is no group of order 3 with all elements involutions.
\end{proof}

\begin{lem}
The graphs $G_9$ and $G_{11}$ are not the principal graphs of subfactors.
\end{lem}
\begin{proof}
This follows immediately from \cite[Corollary 4.7]{1307.5890}.
\end{proof}

\begin{prop}
The graph $G_6$ is not the principal graphs of a subfactor since it does not have a connection.
\end{prop}
\begin{proof}
This follows from \cite[Lemma 3.1]{MR2993924}.
A contradiction follows from the fact that the dimension of the vertex at depth 2 is not equal to the dimension of the non-univalent vertex at depth 4.
\end{proof}

\begin{rem}
$G_7,G_{11}$ are also eliminated by Theorem \ref{thm:3n1}.
\end{rem}

This concludes the classification of subfactors with index exactly 5. We note that all the examples appearing are subgroup subfactors. At the next integer index, 6, much more complicated phenomena arise because 6 is a composite index. 
In particular, we  have many Bisch-Haagerup subfactors. For any normal subgroup $G \vartriangleleft \mathbb Z/2 \mathbb Z * \mathbb Z/3 \mathbb Z$ we have an outer action of the quotient group on the hyperfinite $R$, and the subfactor $R^{\mathbb Z/2 \mathbb Z} \subset R \rtimes \mathbb Z/3\mathbb Z$. However, putting aside the  composite subfactors at  index 6, it remains a possibility that every other example is in fact a subgroup subfactor. This would be initial evidence towards a subfactor analogue of the conjecture that weakly integral fusion categories are weakly group theoretical.



\newcommand{\urlprefix}{}

\bibliographystyle{amsalpha}
{\footnotesize{
\bibliography{../../bibliography/bibliography}
}


\end{document}